\documentclass{amsart}

\usepackage{amssymb}
\usepackage{amsmath}
\usepackage{array}
\usepackage{multirow}

\newcommand{\Lg}{\mbox{$\mathfrak g$}}

\newcommand\fieldsetc{\mathbb}

\newcommand{\R}{\fieldsetc{R}}

\newtheorem{thm}{Theorem}[section]

\theoremstyle{remark}

\newtheorem{quest}{Question}[section]

\newcommand{\D}{\mbox{$\mathcal D$}}
\newcommand{\N}{\mbox{$\mathcal N$}}

\title{Generalized warped products and\\ the $\kappa$-nullity of
Riemannian curvature}
\author{Claudio Gorodski}\thanks{The first named author
  has been partially supported by the grant number 302882/2017-0
  from the \emph{National Council for
  Scientific and Technological Development} (CNPq, Brazil)
  and the project number 2016/23746-6 from the \emph{S\~ao Paulo
  Research Foundation}
  (Fapesp, Brazil).}
\address{Instituto de Matem\'atica e Estat\'\i stica, Universidade de
  S\~ao Paulo, Rua do Mat\~ao, 1010, S\~ao Pulo, SP 05508-090, Brazil.}
\email{gorodski@ime.usp.br}

\author{Felippe Guimar\~aes}\thanks{The second named author
has been supported by the post-doctoral grant 2019/19494-0
from the \emph{S\~ao Paulo Research Foundation}
  (Fapesp, Brazil).}
\address{Instituto de Matem\'atica e Estat\'\i stica, Universidade de
  S\~ao Paulo, Rua do Mat\~ao, 1010, S\~ao Pulo, SP 05508-090, Brazil.}
\email{felippe@impa.br}

\date{\today}

\subjclass[2010]{53C20 (Primary); 53C25, 53C30 (Secondary)}

\begin{document}

\maketitle

{\centering\footnotesize\it Dedicated to Renato de Azevedo Tribuzy 
on the occasion of his 75th anniversary.\par}

\begin{abstract}
In this short survey, we show how two (classes of) known 
examples of inhomogeneous, 
curvature homogeneous Riemannian manifolds with nontrivial $\kappa$-nullity
can be seen as deformations of homogeneous metrics along the vertical
distribution of an integrable Riemannian submersion. We also pose
two open questions.  
\end{abstract}

\section{Introduction}

This is a survey discussing two known examples of complete Riemannian manifolds
with nontrivial $\kappa$-nullity from the stance of generalized
warped products. There are essentially no original results herein, perhaps only
the viewpoint is slightly new.

The idea of nullity
was introduced in case $\kappa=0$
by Chern and Kuiper in~\cite{ChernKuiperNull}, and for
general $\kappa$ by Otsuki~\cite{O}, and later
reformulated and studied by different authors
(see e.g.~\cite{GrayNull,maltzCompl} and, for more recent
work, \cite{DOV} and the references therein). 
Let $M$ be a connected 
Riemannian manifold, and consider the curvature
tensor $R$ of its Levi-Civit\`a connection $\nabla$ with the sign
convention
\[ R(X,Y)Z=\nabla_X\nabla_YZ-\nabla_Y\nabla_XZ-\nabla_{[X,Y]}Z, \]
for vector fields $X$, $Y$, $Z\in\Gamma(TM)$. 
For $\kappa\in\R$, the
\emph{$\kappa$-nullitty distribution} of $M$ is the
variable rank
distribution
$\N_\kappa$ on $M$ defined for each $p\in M$ by
\[ \N_\kappa|_p=\{z\in T_pM:R_p(x,y)z=-\kappa(\langle x,z\rangle_p y-
  \langle y,z\rangle_p x)\quad\mbox{for all $x$, $y\in T_pM$}\}. \]
The number $\nu_\kappa(p):=\dim \N_\kappa|_p$ is called the
\emph{index of $\kappa$-nullity} at~$p$.
We call the orthogonal complement of $\mathcal N_\kappa$ the 
\emph{$\kappa$-conullity distribution} of $M$, and its
dimension at a point~$p\in M$ the \emph{index of $\kappa$-conullity} 
at~$p$, or simply, the $\kappa$-conullity at~$p$.

In case $\kappa=0$ we obtain trivial examples of manifolds
with positive $\nu_0$ simply by taking a Riemannian product
with an Euclidean space, but similar product
examples do not occur if $\kappa\neq0$.  
It is easily seen that
$\nu_\kappa(p)$ is nonzero for at most one value of $\kappa$.
For general~$M$, $\nu_\kappa$ is nonnecessarily constant if nonzero,
but it is a lower semicontinuous function, so there is an open
and dense set of $M$ where $\nu_\kappa$ is locally constant,
and there is
an open subset $\Omega$ of $M$ where $\nu_\kappa$ attains its minimum
value. It is known that $\N_\kappa$ is an autoparallel distribution
on any open set where $\nu_\kappa$ is locally constant and,
in case $M$ is a complete Riemannian manifold,
the leaves of $\N_\kappa$ in $\Omega$ are \emph{complete} totally
geodesic submanifolds of constant curvature $\kappa$~\cite{maltzCompl}.

Having nontrivial $\kappa$-nullity distribution per se 
is usually not a very strong restriction for a Riemannian manifold, 
but, when combined with other conditions (e.g.~homogeneity, irreducibility,
control of the scalar curvature, finite volume, etc.),
 it can led to interesting 
classification/rigidity/non-existence results (see e.g.~\cite{GG}). Further, manifolds
of (maximal) $0$-conullity~$2$ are closely related to semi-symmetric spaces.
Recall that a \emph{semi-symmetric} space is a Riemannian manifold
whose curvature tensor is, at each point, orthogonally equivalent to the
curvature tensor of a symmetric space; the symmetric space may depend on the point.
Conversely, Z. I. Szab\'o~\cite{szabo} showed that every complete irreducible
semi-symmetric space is either locally symmetric or of $0$-conullity~$2$. 
In a related vein,
a Riemannian manifold is called \emph{curvature homogeneous} if 
the curvature tensors at any two points are orthogonally 
equivalent. A Riemannian manifold is said to be \emph{modelled} on 
a given algebraic curvature tensor~$T$ if its curvature tensor is, 
at each point,
orthogonally equivalent to~$T$.  
A curvature homogeneous semi-symmetric space is exactly
a Riemannian manifold whose curvature tensor is modelled
on the curvature tensor of a fixed symmetric space
(see~\cite{BKV} for more details).

\section{A class of generalized warped products and their 
$\kappa$-nullity}\label{sec:warp}

In this section we present formulae for the curvature of a class
of generalized warped products, with a view toward the description 
of two known examples of inhomogeneous manifolds with 
non-trivial $\kappa$-nullity. 
The calculation is straightforward and follows
the lines of~\cite[Ch.~2]{gromoll-walschap}, but it differs
from that book in two respects: herein the function $\varphi$
does not need to be \emph{basic} (that is, constant along the fibers), 
and the Riemannian submersion is assumed
\emph{integrable}, that is, its horizontal distribution is integrable. 

Let $\pi:(M,\langle,\rangle)\to B$ be an integrable Riemannian
submersion.  Denote its Levi-Civit\`a connection and
curvature tensor by $\nabla$, $R$, respectively.
For each
smooth function $\varphi:M\to\R$, we define a new metric
on $M$ by
\begin{equation*}
 \langle u,v\rangle_\varphi=\langle u^h,v^h\rangle+
  e^{2\varphi(p)}\langle u^v,v^v\rangle 
\end{equation*}
for $u$, $v\in T_pM$, where the superscripts denote the horizontal 
and vertical components. Denote the associated Riemannian manifold
by $M_\varphi$, and its Levi-Civit\`a connection and
curvature tensor by $\tilde\nabla$, $\tilde R$, respectively.
Note that $\pi:M_\varphi\to B$ is still an
integrable Riemannian submersion. Denote by $X$, $Y,\ldots$
horizontal vector fields, and by $U$, $V,\ldots$ vertical ones.
The derivation of~$\tilde\nabla$ in terms
of $\nabla$ is a straightforward calculation using the
Koszul formula:
\begin{equation}\label{nabla}
\begin{gathered}
  \tilde\nabla_XY=\nabla_XY,\\
  (\tilde\nabla_VX)^v=(\nabla_VX)^v+X(\varphi)V,\\
  (\tilde\nabla_VX)^h=[V,X]^h\quad \text{($=0$, if $X$ is basic)},\\
  (\tilde \nabla_XV)^v=(\nabla_XV)^v+X(\varphi)V \\
  (\tilde\nabla_XV)^h=0,\\
  (\tilde\nabla_UV)^h=e^{2\varphi}\{(\nabla_UV)^h-\langle U,V\rangle(\nabla\varphi)^h\},\\
 (\tilde \nabla_UV)^v=(\nabla_UV)^v+U(\varphi)V+V(\varphi)U-\langle U,V\rangle(\nabla\varphi)^v,
\end{gathered}
\end{equation}
where $\nabla\varphi$ denotes the gradient of~$\varphi$ with respect
to~$\langle,\rangle$. 

Regarding~$\tilde R$, we need only the following formulae, whose 
derivation, which  uses~(\ref{nabla}), we skip:
\begin{align}\label{curv} \notag
  \tilde R(X,Y)&=R(X,Y)\\ 
  \tilde R(X,V)Y&=R(X,V)Y-X(\varphi)S_YV-Y(\varphi)S_XV\\ \notag
&\qquad\qquad\qquad\qquad+\{\mathrm{Hess}_\varphi(X,Y)+X(\varphi)Y(\varphi)\}V\\ \nonumber
  e^{-2\varphi}\tilde R^h(X,U)V&=R^h(X,U)V+X(\varphi)\sigma(U,V)
  +\langle X,\sigma(U,V)\rangle(\nabla\varphi)^h\\ \notag
                   &\qquad\qquad\qquad\qquad -\langle U,V\rangle\left\{(\nabla_X\nabla\varphi)^h+ X(\varphi)(\nabla \varphi)^h\right\}\\ \nonumber
  \tilde R^v(U,V)X&=R^v(U,V)X-U(\varphi)S_XV+ V(\varphi)S_XU\\ \notag
  &\qquad\qquad\qquad\qquad+\mathrm{Hess}_\varphi(U,X)V-
  \mathrm{Hess}_\varphi(V,X)U.
\end{align}
Here $S$ denotes the shape operators of the fibers, $\sigma$ 
denotes their second fundamental forms, $\mathrm{Hess}_\varphi$
denotes the Hessian of~$\varphi$, and  $(\mathrm{Hess}_\varphi)^h$
is its restriction to the horizontal distribution. 

As an easy consequence of~(\ref{curv}), we deduce the
following result.

\begin{thm}[Nullity of warping]\label{nul-warp}
  Let $\pi:(M,\langle,\rangle)\to B$ be an integrable Riemannian
  submersion and consider $\varphi:M\to\R$.
  \begin{enumerate}
  \item[(i)] If $\pi$ has umbilic fibers (in particular, 
$1$-dimensional fibers) and $\nabla\varphi$
    is vertical, then the horizontal vectors in the $\kappa$-nullity of $M$ 
lie in the $\kappa$-nullity of $M_\varphi$ (here $\kappa$
is arbitrary). 
\item[(ii)] If $S\equiv0$ (so that $\pi$ splits, or $M$ is a local 
product) then 
the horizontal lifts of vectors in the $0$-nullity of $B$,
which in addition lie in  $\ker d\varphi\cap\ker(\mathrm{Hess}_\varphi)^h$, 
belong to the $0$-nullity of $M_\varphi$.
\end{enumerate}
\end{thm}  

\section{Inhomogeneous examples
of metrics with $(-1)$-conullity~$2$}

\setcounter{equation}{0}

In~\cite{TV} Tricerri and Vanhecke quote a problem posed by Gromov, namely,
whether isometry classes of germs of Riemannian metrics 
modelled on the curvature tensor of a given ``irreducible''
homogeneous Riemannian manifold 
depend on a finite number of parameters.  
Schmidt and Wolfson constructed in~\cite{schmidt-wolfson-sl2r}
complete metrics, continuous 
deformations of a certain left-invariant metric on the unimodular 
group $G=SL(2,\R)$, which have constant $(-1)$-conullity~$2$ and 
constant scalar curvature~$-2$; equivalently, these metrics are
curvature homogeneous and are
all modelled on the curvature tensor of the left-invariant
metric on $SL(2,\R)$. Hence 
they answered in the negative Tricerri and Vanhecke's question. 
In this section we discuss those examples 
in light of Theorem~\ref{nul-warp}(i). 

Start with left-invariant vector fields respectively associated to the 
factors of the Iwasawa decomposition $G=NAK$, namely,
\[ Y=\begin{pmatrix}0&1\\0&0\end{pmatrix},\
T=\frac12\begin{pmatrix}1&0\\0&-1\end{pmatrix},\
X=\begin{pmatrix}0&1\\-1&0\end{pmatrix}. \]
Then we have the bracket relations
\begin{equation}\label{sl2r}
  [X,Y]=2T,\ [T,X]=-X+2Y,\ [T,Y]=Y.
  \end{equation}
Note that $Y$ and $T$ span an involutive distribution. 
Consider the left-invariant metric obtained by declaring 
the frame~$X$, $Y$, $T$ orthonormal. Then one easily computes that 
\begin{equation}\label{lc0}
  \begin{gathered}
    \nabla_TT=\nabla_TX=\nabla_TY=0,\ \nabla_XT=X-2Y,\ \nabla_YT=-Y,\\
 \nabla_XX=-T,\ \nabla_YY=T,\ \nabla_XY=2T,\
  \nabla_YX=0.
\end{gathered}
\end{equation}
We deduce that the distribution $\D$ spanned by $Y$ and $T$ is autoparallel,
its leaves are isometric to real hyperbolic planes $\R H^2(-1)$
of curvature $-1$, and $G$ has $(-1)$-nullity distribution of rank~$1$,
spanned by $T$. We note that the sectional curvature of the plane 
spanned by $X$, $Y$ is $1$,
which implies constancy of the scalar curvature of $G$, equal to $-2$. 

Consider now the projection $\pi:G\to G/K=NA=\R H^2(-1)$. 
The vector field $X$ is Killing transversal, that is, its flow 
is isometric on $\D$ (this follows from~(\ref{lc0})), so its flow lines 
induce a Riemannian foliation of $G$ and $\pi$ is the projection onto the 
space of leaves. In particular, $\pi$ is a Riemannian submersion,
the horizontal distribution being spanned by $Y$, $T$ and thus integrable, 
and we can apply Theorem~\ref{nul-warp}(i). For every smooth
$\varphi:G\to\R$ with a vertical gradient, $G_\varphi$ also 
has $(-1)$-nullity distribution of rank~$1$,
spanned by $T$. It is easily seen that the scalar curvature of 
$G_\varphi$ is likewise~$-2$ (this also follows from the last 
equation in~(\ref{eqns})). Since $\nabla\varphi$ is vertical,
$\varphi$ is actually a function on the compact group $K$ and hence bounded. 
This easily implies that every divergent curve in $G_\varphi$ has 
infinite length and hence $G_\varphi$ is complete as a Riemannian manifold. 
We remark that inhomogeneous Schmidt-Wolfson metrics do not cover a manifold
of finite volume~\cite[Thm.~1.1]{sw}.  

How general are these examples? Namely, 
assume $M$ is a given Riemannian $3$-manifold with minimal
$(-1)$-nullity~$1$
and constant scalar curvature. 
Around any point of minimal nullity, 
we can find a local orthonormal frame $T$, $X$, $Y$, where 
$T$ spans the nullity distribution, satisfying (we refer to~\cite{GG}
for this construction, which is not difficult):
\begin{equation}\label{lc}
  \begin{gathered}
    \nabla_TT=\nabla_TX=\nabla_TY=0,\ \nabla_XT=X-2FY,\ \nabla_YT=-Y,\\
 \nabla_XX=-T+\alpha Y,\ \nabla_YY=T+\beta X,\ \nabla_XY=2FT-\alpha X,\
  \nabla_YX=-\beta Y.
\end{gathered}
\end{equation}
for some locally defined smooth functions $\alpha$, $\beta$. 
The bracket relations follow:
\begin{equation*}
  [X,Y]=2FT-\alpha X+\beta Y,\ [T,X]=-X+2FY,\ [T,Y]=Y.
  \end{equation*}
Next, the curvature relations
\[ \langle R(X,Y)X,Y\rangle=-k, \]
where $k$ is the sectional curvature of the
plane spanned by $X$, $Y$, and 
\[  \langle R(X,Y)X,T\rangle=\langle R(T,Y)X,Y\rangle=
  \langle R(X,Y)Y,T\rangle=0, \]
yield the equations
\begin{equation}\label{eqns}
\begin{aligned}
  \alpha &= -\beta F,\\ 
  T(\beta)&=\beta,\\
  Y(F)&=-\beta(1+F^2),\\
  X(\beta)-FY(\beta) &= k-1.
\end{aligned}                      
\end{equation}

Note that $Y$ and $T$ span an involutive distribution, with 
integral manifolds isometric to~$\R H^2(-1)$. Further, in view
  of~(\ref{lc}) these
  hyperbolic planes are totally geodesic if and only if $\beta=0$
if and only if the foliation by integral curves of $X$ is 
Riemannian (indeed a \emph{polar foliation}, see e.g.~\cite{thorbergsson}) 
if and only if the projection of $SL(2,\R)$ onto the 
space of leaves is a Riemannian submersion; 
  equivalently, $Y(F)=0$. \emph{Assume this is the case};
  we work on the open set $U$ defined by $F\neq0$, where
  we can take a new frame $\tilde T=T$, $\tilde X=F^{-1}X$, $\tilde Y=Y$.
  The new frame satisfies 
  \[    [\tilde X,\tilde Y]=2\tilde T,\ [\tilde T,\tilde X]=-\tilde X+2\tilde Y,\ [\tilde T,\tilde Y]=\tilde Y. \]
Comparing with~(\ref{sl2r}), 
  by Lie's third theorem, the metric defined on $U$ that
  declares $\tilde X$, $\tilde Y$, $\tilde T$ orthonormal is locally isometric
  to a left-invariant metric on $SL(2,\R)$. We deduce that the
  original metric on $U$ is locally isometric to a Schmidt-Wolfson metric
on $SL(2,\R)$, that is, a deformation as above.  

\begin{quest}
Are there exist complete solutions to~(\ref{eqns}) with $\beta\neq0$?
\end{quest}  

\begin{quest}
Are there other examples of polar foliations of dimension one 
(or with umbilic leaves) on non-symmetric spaces, other than $SL(2,\R)$?
\end{quest}

\section{Inhomogeneous examples of metrics with $0$-conullity~$2$}
\setcounter{equation}{0}

In 1968, Nomizu conjectured that every complete irreducible semi-symmetric 
space  of dimension greater than or equal to three is locally symmetric.
His conjecture was refuted by Takagi and Sekigawa in 1972, who constructed
counterexamples.  Here we are concerned with another two of
Sekigawa's counterexamples~\cite{sekigawa}, later generalized  
by Kowalski, Tricerri and Vanhecke~\cite{KTV} (these were later 
further generalized, see~\cite{BKV} for the full range of generalizations).

We start with an almost Abelian Lie group 
(that is, a Lie group with a codimension one ideal) 
$G=S^1\ltimes_\rho V$, 
where the representation $\rho:S^1\to SO(V)$ is orthogonal. 
Its Lie algebra is $\Lg=\R\ltimes_{d\rho}V$. We fix 
$\xi\in\R$ and consider $d\rho(\xi)\in\mathfrak{so}(V)$. 
We also equip $G$ with a left-invariant Riemannian metric
that makes $\xi$ unit, and $\xi$ and $V$ orthogonal. 
Koszul's formula for the 
Levi-Civit\`a connection immediately 
yields:
\begin{equation}\label{aa-lc}
 \nabla_\xi\xi=\nabla_X\xi=\nabla_XY=0,\ \nabla_\xi X=d\rho(\xi)X, 
\end{equation}
for all~$X$, $Y\in V$. It immediately follows that the metric is flat.

Next we are going to deform the metric on $G$, as
in Theorem~\ref{nul-warp}(ii). 
The projection $\pi:G\to V$ is obviously an integrable Riemannian 
submersion. We need to construct the function $\varphi$.   
Let $Z\in V$ be a fixed unit vector, and consider the associated 
height function $h(X)=\langle X,Z\rangle$ for $X\in V$. Extend~$h$ to a 
smooth function $t:G\to\R$ by setting $t(g)=h(\rho(e^{-i\theta})Y)$,
for $g=(e^{i\theta},Y)\in G$. Below we are going to use the 
following calculation: for any $X\in V$, we have
\[ g\exp(sX)=(e^{i\theta},Y)(1,sX)=(e^{i\theta},s\rho(e^{i\theta})X+Y), \]
so
\begin{align*}
  X_g(t)&=\frac{d}{ds}\Big|_{s=0}t(g\exp(sX))\\
        &=\frac{d}{ds}\Big|_{s=0}\langle sX+\rho(e^{-i\theta})Y,Z\rangle\\
        &=\langle X,Z\rangle,
\end{align*}
and hence $(\nabla t)^h=Z$. Finally, define $\varphi=\log f$, where 
\[  f= c_1 e^{at}+ c_2 e^{-at} \]
for constants $c_1$, $c_2\geq0$ and $a>0$. 

Back to Theorem~\ref{nul-warp}(ii): of course $\pi$ splits, and the 
$0$-nullity of the base $V$ is $V$ itself. We have
\[ (\nabla\varphi)^h= \varphi_t(\nabla t)^h=\varphi_tZ, \]
so the horizontal part of $\ker d\varphi$ is the left-invariant, 
horizontal hyperplane distribution $Z^\perp$. 
For $X\in Z^\perp$ and $Y\in V$, we have
\begin{align*}
 \mathrm{Hess}_\varphi(Y,X)&=YX(\varphi)-\nabla_YX(\varphi)\\
  &=YX(\varphi)\quad\text{(since $\nabla_YX=0$)} \\
  & = 0 \quad\text{(since $X(\varphi)=0$),} 
\end{align*}
that is, $(\ker\mathrm{Hess}_\varphi)^h\supset\ker d\varphi$. 
It follows from the theorem that the $0$-nullity distribution of $G_\varphi$ 
is $\N_0=Z^\perp$. In particular, $G_\varphi$ has constant $0$-conullity~$2$. 

Consider the unit vector field 
$\tilde\xi=\frac1f\xi$ in the metric $\langle,\rangle_\varphi$. 
From formulae~(\ref{curv}), we obtain:
\begin{equation}\label{curv2}
\begin{gathered}
\tilde R(X,Y)=0,\\
\tilde R(\tilde\xi,X)Y=-a^2\langle X,Z\rangle\langle Y,Z\rangle\tilde\xi,\\
\tilde R(\tilde\xi,X)\tilde\xi=a^2\langle X,Z\rangle Z,
\end{gathered}
\end{equation}
for $X$, $Y\in V$. From here we see that the sectional curvature
$\tilde K(\tilde\xi,Z)=-a^2$, and hence $G_\varphi$ has constant
(negative) scalar curvature $-2a^2$. It follows that $G_\varphi$ is curvature
homogeneous, semi-symmetric, and modelled on $\R H^2(-a^2)\times\R^{m-1}$,
where $m=\dim V$. 

Next we prove completeness of $G_\varphi$.
A Riemannian manifold is complete if and only if every divergent 
curve has infinite length. Let $\gamma$ be a divergent curve in $G_\varphi$. 
Since $S^1$ is compact, also $\pi\circ\gamma$ is divergent,
so by completeness of 
the base $V$ we have that the length
$L(\gamma)\geq L(\pi\circ\gamma)=\infty$. 
This proves that $G_\varphi$ is complete. 

Finally, we state a sufficient condition for $G_\varphi$ to be 
irreducible as a Riemannian manifold. Let $\hat V=\R\oplus V$ denote the 
tangent space to $G_\varphi$ at the basepoint. 
We are going to prove that 
if $Z$ is a cyclic vector for the representation $d\rho$, then 
the identity component of the holonomy group of $G_\varphi$ 
at the basepoint is the full rotation group $SO(\hat V)$; 
in particular, this holonomy group 
acts irreducibly on $\hat V$. 

It is convenient to use the standard identification 
$\bigwedge^2\hat V\cong\mathfrak{so}(\hat V)$ given by 
$u\wedge v\mapsto\langle u,\cdot\rangle_\varphi v-
\langle v,\cdot\rangle_\varphi u$. Then the Lie bracket
$[v\wedge w,u\wedge v]=u\wedge w$, for $u$, $v$, $w\in\hat V$. 
Recall that $m=\dim V$, let $T=d\rho(\tilde\xi)$, and 
denote by $\hat V_k$ the span of $\tilde\xi,Z,\ldots,T^kZ$,
where $k=0,\ldots,m-1$. 

It is enough to prove that, under the 
assumption that $Z$ is a cyclic vector, the
Lie algebra of the 
infinitesimal holonomy at the basepoint is $\mathfrak{so}(\hat V)$;
indeed we will show by induction on~$k_0$ that the Lie subalgebra
generated by 
\[ \{\tilde\nabla_{\tilde\xi}^k\tilde R(\tilde\xi,Z)\in\mathfrak{so}(\hat V):k=0,\ldots,k_0\} \] 
contains $\mathfrak{so}(\hat V_{k_0})$, where $k_0=0,\ldots,m-1$.   
In turn, this result will follow from the following claims,
where a decomposable element of $\bigwedge^2\hat V$ is said to have 
order $<k$ if it only involves $\tilde\xi$, $Z,\ldots,T^{k-1}Z$:

\noindent\textsc{Claim 1.} 
$\tilde\nabla_{\tilde\xi}^k\tilde R(\tilde\xi,Z)=a^2(\tilde\xi\wedge T^kZ)+
\mbox{terms of order $<k$.}$

\noindent\textsc{Claim 2.} 
$\tilde\nabla_{\tilde\xi}^k\tilde R(\tilde\xi,X)=\mbox{sum of terms of order $<k$}$,
where $X\in Z^\perp$. 

\noindent\textsc{Claim 3.} 
$\tilde\nabla_{\tilde\xi}^k\tilde R(X,Y)=\mbox{sum of terms of order $<k$}$,
where $X$, $Y\in V$. 

We combine~(\ref{nabla}) with~(\ref{aa-lc}) to write (by the way,
a left-invariant $X\in V$ is generally not basic with respect to
$\pi$):
\begin{equation}\label{nabla-tilde}
\begin{gathered}
\tilde\nabla_XY=0,\\
\tilde\nabla_{\tilde\xi}\tilde\xi=-\varphi_tZ,\\
\tilde\nabla_{\tilde\xi}X=TX+\varphi_t\langle X,Z\rangle\tilde\xi,
\end{gathered}
\end{equation}
where $X$, $Y\in V$. 
Also, from~(\ref{curv2}) we see that 
\begin{equation*}
\begin{gathered}
 \tilde R(\tilde\xi,Z)=a^2\tilde\xi\wedge Z,\\
 \tilde R(\tilde\xi,X)=0,\quad\text{if $X\in Z^\perp$,}\\ 
\tilde R(X,Y)=0,\quad\text{for $X$, $Y\in V$.} 
\end{gathered}
\end{equation*}
This already gives the initial case $k=0$ in claims~$1$ through~$3$;
their proof now follows from a straightforward 
induction on~$k$, using~(\ref{nabla-tilde}) in 
\begin{align*}
\tilde\nabla^{k+1}_{\tilde\xi}\tilde R(\tilde \xi,Z)&=\tilde\nabla_{\tilde\xi}
(\tilde\nabla_{\tilde \xi}^k\tilde R)(\tilde\xi,Z)\\
&=\tilde\nabla_{\tilde\xi}(\tilde\nabla_{\tilde \xi}^k\tilde R(\tilde\xi,Z))
-\tilde\nabla^k_{\tilde\xi}\tilde R(\tilde\nabla_{\tilde\xi}\tilde\xi,Z)
-\tilde\nabla^k_{\tilde\xi}\tilde R(\tilde\xi,\tilde\nabla_{\tilde\xi}Z),
\end{align*}
and in similar identities for 
$\tilde\nabla^{k+1}_{\tilde\xi}\tilde R(\tilde \xi,X)$ (where~$X\in Z^\perp$) and
$\tilde\nabla^{k+1}_{\tilde\xi}\tilde R(X,Y)$ (where~$X$, $Y\in V$). 

This completes the proof of irreducibility 
of $G_\varphi$ under the assumption that $Z$ is a cyclic vector,
and the construction of the examples. We finish with two remarks: first, it is 
elementary to see that cyclic vectors
exist in $V$ if and only if $\rho$ is multiplicity free (that is, 
it has pairwise different weights as a representation of $S^1$),
and in such a case, the cyclic vectors are exacly those vectors
with nonzero components in all irreducible components; and 
second, in~\cite{KTV}
it was further shown that cyclicity of the vector $Z$ is also a 
necessary condition for the irreducibility of~$G_\varphi$. 

\providecommand{\bysame}{\leavevmode\hbox to3em{\hrulefill}\thinspace}
\providecommand{\MR}{\relax\ifhmode\unskip\space\fi MR }
\providecommand{\MRhref}[2]{%
  \href{http://www.ams.org/mathscinet-getitem?mr=#1}{#2}
}
\providecommand{\href}[2]{#2}

\end{document}